\documentclass[11pt]{article}
\usepackage{}
\usepackage{mathrsfs}
\usepackage{amsfonts}
\usepackage{amsfonts,amssymb,mathrsfs}
\usepackage{amsmath,amscd}
\usepackage[dvips]{graphicx}
\usepackage[all]{xy}
\usepackage{graphicx}
\usepackage{fancyhdr}
\usepackage{mathptm,pslatex}
\usepackage{amsthm}

\begin{document}

\sloppy
\newtheorem{Def}{Definition}[section]
\newtheorem{Bsp}{Example}[section]
\newtheorem{Prop}[Def]{Proposition}
\newtheorem{Theo}[Def]{Theorem}
\newtheorem{Lem}[Def]{Lemma}
\newtheorem{Koro}[Def]{Corollary}
\theoremstyle{definition}
\newtheorem{Rem}[Def]{Remark}

\newcommand{\add}{{\rm add}}
\newcommand{\gd}{{\rm gl.dim }}
\newcommand{\dm}{{\rm dom.dim }}
\newcommand{\E}{{\rm E}}
\newcommand{\Mor}{{\rm Morph}}
\newcommand{\End}{{\rm End}}
\newcommand{\ind}{{\rm ind}}
\newcommand{\rsd}{{\rm res.dim}}
\newcommand{\rd} {{\rm rep.dim}}
\newcommand{\ol}{\overline}
\newcommand{\rad}{{\rm rad}}
\newcommand{\soc}{{\rm soc}}
\newcommand{\su}{{\rm sup}}
\renewcommand{\top}{{\rm top}}
\newcommand{\pd}{{\rm proj.dim}}
\newcommand{\id}{{\rm inj.dim}}
\newcommand{\Fac}{{\rm Fac}}
\newcommand{\fd} {{\rm fin.dim }}
\newcommand{\DTr}{{\rm DTr}}
\newcommand{\cpx}[1]{#1^{\bullet}}
\newcommand{\D}[1]{{\mathscr D}(#1)}
\newcommand{\Dz}[1]{{\mathscr D}^+(#1)}
\newcommand{\Df}[1]{{\mathscr D}^-(#1)}
\newcommand{\Db}[1]{{\mathscr D}^b(#1)}
\newcommand{\C}[1]{{\mathscr C}(#1)}
\newcommand{\Cz}[1]{{\mathscr C}^+(#1)}
\newcommand{\Cf}[1]{{\mathscr C}^-(#1)}
\newcommand{\Cb}[1]{{\mathscr C}^b(#1)}
\newcommand{\K}[1]{{\mathscr K}(#1)}
\newcommand{\Kz}[1]{{\mathscr K}^+(#1)}
\newcommand{\Kf}[1]{{\mathscr  K}^-(#1)}
\newcommand{\Kb}[1]{{\mathscr K}^b(#1)}
\newcommand{\modcat}{\ensuremath{\mbox{{\rm -mod}}}}
\newcommand{\Modcat}{\ensuremath{\mbox{{\rm -Mod}}}}
\newcommand{\stmodcat}[1]{#1\mbox{{\rm -{\underline{mod}}}}}
\newcommand{\pmodcat}[1]{#1\mbox{{\rm -proj}}}
\newcommand{\Pmodcat}[1]{#1\mbox{{\rm -Proj}}}
\newcommand{\imodcat}[1]{#1\mbox{{\rm -inj}}}
\newcommand{\opp}{^{\rm op}}
\newcommand{\otimesL}{\otimes^{\rm\bf L}}
\newcommand{\rHom}{{\rm\bf R}{\rm Hom}}
\newcommand{\projdim}{\pd}
\newcommand{\Hom}{{\rm Hom}}
\newcommand{\Coker}{{\rm coker}\,\,}
\newcommand{ \Ker  }{{\rm Ker}\,\,}
\newcommand{ \Img  }{{\rm Im}\,\,}
\newcommand{\Ext}{{\rm Ext}}
\newcommand{\StHom}{{\rm \underline{Hom} \, }}

\newcommand{\gm}{{\rm _{\Gamma_M}}}
\newcommand{\gmr}{{\rm _{\Gamma_M^R}}}

\def\vez{\varepsilon}\def\bz{\bigoplus}  \def\sz {\oplus}
\def\epa{\xrightarrow} \def\inja{\hookrightarrow}

\newcommand{\lra}{\longrightarrow}
\newcommand{\lraf}[1]{\stackrel{#1}{\lra}}
\newcommand{\ra}{\rightarrow}
\newcommand{\dk}{{\rm dim_{_{k}}}}

{\Large \bf
\begin{center}
Generalized Auslander-Reiten conjecture and derived equivalences
\end{center}}
\medskip

\centerline{{\bf  Shengyong Pan}}
\begin{center} Department of Mathematics,\\
 Beijing Jiaotong University,  Beijing 100044,\\
People's Republic of China\\ E-mail:shypan@bjtu.edu.cn \\
\end{center}
\bigskip

\renewcommand{\thefootnote}{\alph{footnote}}
\setcounter{footnote}{-1} 
\renewcommand{\thefootnote}{\alph{footnote}}
\setcounter{footnote}{-1} 


\begin{abstract}
In this note, we prove that the generalized Auslander-Reiten
conjecture is preserved under derived equivalences between Artin
algebras.
\end{abstract}

\section{Introduction}

In the representation theory of Artin algebras, one of the most
important open problems is the Nakayama conjecture which predicts
that an Artin algebra $A$ is self-injective provided that all terms
in a minimal injective resolution of $A$ are projective. Mueller
\cite{M} in the late sixties proved that the Nakayama conjecture
holds for every Artin algebra if and only if for any Artin algebra
$A$ any finitely generated generator-cogenerator $M$, the vanishing
$\Ext^{n}_A(M,M)=0$, for $n\geq1$, implies that $M$ is projective.
In this connection and based on Mueller's result, Auslander-Reiten
proposed several stronger conjectures, and in particular the
following:

(ARC) \quad Let $M$ be a finitely generated module over an Artin
algebra $A$ such that $\Ext_A^i(M,M)=0=\Ext_A^i(M,A), $ for
$i\geq1$. Then $M$ is projective.

The above conjecture, widely known as {\it the Auslander-Reiten
Conjecture}, implies the Nakayama conjecture, and, as
Auslander-Reiten proved, it holds in all cases where the Nakayama
conjecture is known to be true. It should be noted that the above
conjectures, which are trivial consequences of the finitistic
dimension conjecture, are still open. Auslander-Reiten conjecture
has been verified for some special classes of Artin algebras and
commutative Noetherian rings \cite{CH,W1,W3}.

In this note, we consider the following generalization version of
Auslander and Reiten conjecture which can be stated as follows:

(GARC) \quad Let $A$ be an Artin algebra. Let $X$ be a finitely
generated $A$-module and $r$ a non-negative integer. If
$$
\Ext_A^i(M,M)=0=\Ext_A^i(M,A),
$$
for $i> r$, then $\pd(X)\leq r$, where $\pd(X)$ is the projective
dimension of $X$.

In case $r=0$, (GARC) is (ARC). In \cite {W1}, the generalized
Auslander-Reiten conjecture holds for Artin algebras for which any
finitely generated module has an ultimately closed projective
resolution. It also holds for all algebras which satisfy the
Auslander-Reiten conjecture. In \cite{W2}, it was proved that the
generalized Auslander-Reiten conjecture is stable under tilting
equivalences.

The aim of this note is to show that the generalized Auslander
Reiten conjecture is preserved by derived equivalences, as it was
done for the finiteness of finitistic dimension conjecture
\cite{PX}. Then we generalize the main result of \cite{W2},
answering in the affirmative a question of Wei.

Our main result reads as follows:

\begin{Theo}
Suppose that $A$ and $B$ are Artin algebras. Assume that $A$ and $B$
are derived equivalent. Then $A$ satisfies the generalized
Auslander-Reiten conjecture if and only if so does $B$. \label{thm1}
\end{Theo}

In view of the importance of the Nakayama conjecture and the
(Generalized) Auslander-Reiten conjecture, it is highly desirable to
have as much as possible information about classes of algebras
satisfying the conjectures. The main result indicates that the
validity of the generalized Auslander-Reiten conjecture for an Artin
algebra depends on its derived equivalence class and in this way one
produces further classes of algebras satisfying the conjecture.

The contents of this paper are organized as follows. In Section
\ref{pre}, we recall some definitions and notations on derived
categories and derived equivalences. In Section \ref{i}, we prove
our main result, Theorem \ref{thm1}.\\

\noindent{\bf Acknowledgements.} The author would like to
acknowledge the Fundamental Research Funds for the Center University
(2011JBM131) and postdoctoral granted financial support from China
Postdoctoral Science Foundation (20100480188), during which this
work was carried out. The author also would like to thank the
referee for his/her helpful comments that improve the paper.

\section{Preliminaries}\label{pre}
In this section, we shall recall some definitions and notations on
derived categories and derived equivalences, and basic results which
are needed in the proofs of our main results.

Let $\mathscr{A}$ be an abelian category. For two morphisms $\alpha:
X\ra Y$ and $\beta: Y\ra Z$, their composition is denoted by
$\alpha\beta$. An object $X\in\mathscr{A}$ is called a additive
generator for $\mathscr{A}$ if $\add(X)=\mathscr{A}$, where
$\add(X)$ is the additive subcategory of $\mathscr{A}$ consisting of
all direct summands of finite direct sums of the copies of $X$. A
complex $\cpx{X}=(X^i,d_{X}^i)$ over $\mathscr{A}$ is a sequence of
objects $X^i$ and morphisms $d_{X}^i$ in $\mathscr{A}$ of the form:
$\cdots \ra X^i\stackrel{d^i}\ra X^{i+1}\stackrel{d^{i+1}}\ra
X^{i+2}\ra\cdots$, such that $d^id^{i+1}=0$ for all
$i\in\mathbb{Z}$. If $\cpx{X}=(X^i,d_{X}^i)$ and
$\cpx{Y}=(Y^i,d_{Y}^i)$ are two complexes, then a morphism $\cpx{f}:
\cpx{X}\ra\cpx{Y}$ is a sequence of morphisms $f^i: X^i\ra Y^i$ of
$\mathscr{A}$ such that $d^i_{X}f^{i+1}=f^id^i_{Y}$ for all
$i\in\mathbb{Z}$. The map $\cpx{f}$ is called a chain map between
$\cpx{X}$ and $\cpx{Y}$. The category of complexes over
$\mathscr{A}$ with chain maps is denoted by $\C{\mathscr{A}}$. The
homotopy category of complexes over $\mathscr{A}$ is denoted by
$\K{\mathscr{A}}$ and the derived category of complexes is denoted
by $\D{\mathscr{A}}$.

Let $R$ be a commutative Artin ring, and let $A$ be an Artin
$R$-algebra. We denote by $A$-Mod and $A$-mod the categories of left
$A$-modules and finitely generated left $A$-modules, respectively.
The full subcategories of $A$-Mod and $A$-mod consisting of
projective modules and finitely generated projective modules are
denoted by $A$-Proj and $_A\mathcal {P}$, respectively. Denote by
$_{A}\mathcal {X}^{>r}$ the category of $A$-modules $X$ satisfied
$\Ext_{A}^{i}(X,A)=0$ for $i> r$, where $r$ is a non-negative
integer. Particularly, if an $A$-module $X$ satisfies
$\Ext_{A}^{i}(X,A)=0$ for $i>0$, then $X$ is said to be a
Cohen-Macaulay $A$-module.

Recall that a homomorphism $f: X\ra Y$ of $A$-modules is called a
radical map provided that for any $A$-module $Z$ and homomorphisms
$g: Y\ra Z$ and $h: Z\ra X$, the composition $hfg$ is not an
isomorphism. A complex of $A$-modules is called a radical complex if
its differential maps are radical maps. Let $\Kf{A}$ and $\Kb{A}$
denote the homotopy category of bounded above and bounded complexes
of $A$-modules, respectively. We denote by $\Df{A}$ and $\Db{A}$ the
derived category of bounded above and bounded complexes of
$A$-modules, respectively.

The fundamental theory on derived equivalences has been established
by Rickard \cite{Ri1}.

\begin{Theo}$\rm \cite[Therem\;6.4]{Ri1}$\label{2.1}
Let $A$ and $B$ be rings. The following conditions are equivalent.

$(i)$ $\Db{A\Modcat}$ and $\Db{B\Modcat}$ are equivalent as
triangulated categories.

$(ii)$ $\Kf{\Pmodcat{A}}$ and $\Kf{\Pmodcat{B}}$ are equivalent as
triangulated categories.

$(iii)$ $\Kb{\Pmodcat{A}}$ and $\Kb{\Pmodcat{B}}$ are equivalent as
triangulated categories.

$(iv)$ $\Kb{_{A}\mathcal {P}}$ and $\Kb{_{B}\mathcal {P}}$ are
equivalent as triangulated categories.

$(v)$ $B$ is isomorphic to $\End_{\Db{A}}(\cpx{T})$ for some complex
$\cpx{T}$ in $\Kb{_{A}\mathcal {P}}$ satisfying

         \qquad $(1)$ $\Hom_{\Db{A}}(\cpx{T},\cpx{T}[n])=0$
         for all $n\neq 0$.

         \qquad $(2)$ $\add(\cpx{T})$, the category of direct summands of
          finite direct sums of copies of $\cpx{T}$, generates
          $\Kb{_{A}\mathcal {P}}$ as a triangulated category.
\end{Theo}

\noindent{\bf Remarks.} (1) The rings $A$ and $B$ are said to be
derived equivalent if $A$ and $B$ satisfy the conditions of the
above theorem. The complex $\cpx{T}$ in Theorem 2.1 is called a
tilting complex for $A$.

(2) By \cite[Corollary 8.3]{Ri1}, two Artin $R$-algebras $A$ and $B$
are said to be  derived equivalent if their derived categories
$\Db{A}$ and $\Db{B}$ are equivalent as triangulated categories. By
Theorem 2.1, Artin algebras $A$ and $B$ are derived equivalent if
and only if $B$ is isomorphic to the endomorphism algebra of a
tilting complex $\cpx{T}$. If $\cpx{T}$ is a tilting complex for
$A$, then there is an equivalence $F: \Db{A}\ra\Db{B}$ that sends
$\cpx{T}$ to $B$. On the other hand, for each derived equivalence
$F: \Db{A}\ra\Db{B}$, there is an associated tilting complex
$\cpx{T}$ for $A$ such that $F(\cpx{T})$ is isomorphic to $B$ in
$\Db{B}$.

\section{Generalized Auslander-Reiten conjecture is invariant under derived equivalences}\label{i}
In this section, we shall prove Theorem \ref{thm1}.

We shall give the proof as a sequence of lemmas.

Now we suppose that $A$ and $B$ are Artin algebras. Let $F:
\Db{A}\lra \Db{B}$ be a derived equivalence and let $\cpx{P}$ be the
tilting complex associated to $F$. Without loss of generality, we
assume that $\cpx{P}$ is a radical complex of the following form
$$
0\ra P^{-n}\ra P^{-n+1} \ra \cdots\ra P^{-1}\ra P^{0}\ra 0.
$$

Then we have the following fact.

\begin{Lem} $\rm\cite[lemma\, 2.1]{hx}$
Let $F: \Db{A}\lra \Db{B}$ be a derived equivalence between Artin
algebras $A$ and $B$. Then there is a tilting complex
$\bar{P}^{\bullet}$ for $B$ associated to the quasi-inverse of $F$
of the form
$$
0\ra \bar{P}^{0}\ra \bar{P}^{1} \ra \cdots\ra \bar{P}^{n-1}\ra
\bar{P}^{n}\ra 0,
$$ with the differential being radical maps.
\end{Lem}

Suppose that $\cpx{X}$ is a complex of $A$-modules. We define the
following truncations:

$\tau_{\geq 1}(\cpx{X}): \cdots\ra0\ra0\ra X^{1}\ra X^{2}\ra\cdots$,

$\tau_{\leq 0}(\cpx{X}): \cdots\ra X^{-1}\ra X^{0}\ra 0\ra0\cdots$.

The following lemma, proved in \cite[Lemma\,2.1]{PX}, will be used
frequently in our proofs below.

\begin{Lem}\label{PX} Let $m,t,d \in \mathbb{N}$, $X^{\bullet}, Y^{\bullet} \in
\Kb{A}$. Assume that $X^{i}=0$ for $i<m$, $Y^{j}=0$ for $j>t$, and
$\Ext^{l}(X^{i},Y^{j})=0$ for all $i,j\in \mathbb{N}$ and $l\geq d$.
Then $\Hom_{\Db{A}}(X^{\bullet},Y^{\bullet}[l])=0$ for $l\geq
d+t-m$.
\end{Lem}

The following lemma is inspired by \cite[lemma 3.3]{P}, we have a
variation.

\begin{Lem} Let $F:
\Db{A}\lra \Db{B}$ be a derived equivalence between Artin algebras
$A$ and $B$, and let $G$ be the quasi-inverse of $F$. Suppose that
$\cpx{P}$ and $\bar{P}^{\bullet}$ are the tilting complexes
associated to $F$ and $G$, respectively. Let $r$ be a non-negative
integer. Then

$(i)$ For $X\in _A\mathcal {X}^{>r}$, the complex $F(X)$ is
isomorphic in $\Db{B}$ to a radical complex $\bar{P}^{\bullet}_{X}$
of the form
$$
0\ra \bar{P}_{X}^{0}\ra \bar{P}_{X}^{1} \ra \cdots\ra
\bar{P}_{X}^{n-1}\ra \bar{P}_{X}^{n}\ra 0
$$
with $\bar{P}_{X}^{0}\in_B\mathcal {X}^{>r}$ and $\bar{P}_{X}^{i}$
projective $B$-modules for $1\leq i\leq n$.

$(ii)$ For $Y\in_B\mathcal {X}^{>r}$, the complex $G(Y)$ is
isomorphic in $\Db{A}$ to a radical complex $P^{\bullet}_{Y}$ of the
form

$$
0\ra P_{Y}^{-n}\ra P_{Y}^{-n+1} \ra \cdots\ra P_{Y}^{-1}\ra
P_{Y}^{0}\ra 0
$$
with $P_{Y}^{-n}\in_A\mathcal {X}^{>r}$ and $P_{Y}^{i}$ projective
$A$-modules for $-n+1\leq i\leq 0$.\label{2.3}
\end{Lem}

{\bf Proof}. We only to show the first case. The proof of ($ii$) is
similar to that of ($i$).

($i$) For an $A$-module $X$ with $\Ext_A^i(X,A)=0$ for $i>r$, by
\cite[Lemma 3.1]{hx}, we see that the complex $F(X)$ is isomorphic
in $\Db{B}$ to a complex $\bar{P}^{\bullet}_{X}$ of the form
$$
0\ra \bar{P}_{X}^{0}\ra \bar{P}_{X}^{1} \ra \cdots\ra
\bar{P}_{X}^{n-1}\ra \bar{P}_{X}^{n}\ra 0,
$$
with $\bar{P}_{X}^{i}$ projective $B$-modules for $i>0$. We only
need to show that $\Ext_B^{i>r}(\bar{P}_{X}^{0},B)=0$, that is,
$\Ext^{i}_{B}(\bar{P}_{X}^{0},B)=0$ for $i> r$, where $r$ is a
non-zero integer. Indeed, there exists a distinguished triangle
$$\bar{P}^{+}_{X}\ra\bar{P}^{\bullet}_{X}\ra\bar{P}^{0}_{X}\ra
\bar{P}^{+}_{X}[1]$$ in $\Kb{B}$, where $\bar{P}^{+}_{X}$ denotes
the complex $\tau_{\geq 1}(\cpx{\bar{P}_{X}})$. For each
$i\in\mathbb{Z}$, applying the functor $\Hom_{\Db{B}}(-,B[i])$ to
the above distinguished triangle, we get an exact sequence
\begin{eqnarray*}
\cdots\ra\Hom_{\Db{B}}(\bar{P}^{+}_{X}[1],B[i])\ra
\Hom_{\Db{B}}(\bar{P}^{0}_{X},B[i])\ra
\Hom_{\Db{B}}(\bar{P}^{\bullet}_{X},B[i])\\\ra
\Hom_{\Db{B}}(\bar{P}^{+}_{X},B[i])\ra\cdots. \end{eqnarray*} On the
other hand, $\Hom_{\Db{B}}(\bar{P}^{+}_{X},B[i])\simeq
\Hom_{\Kb{B}}(\bar{P}^{+}_{X},B[i])=0$ for $i> r$. By Lemma \ref{PX}
and $\End^{i}_{A}(X,A)=0$ for $i>r$, we get
$\Hom_{\Db{B}}(\bar{P}^{\bullet}_{X},B[i])\simeq\Hom_{\Db{A}}(X,P^{\bullet}[i])=0$
for all $i> r$. Consequently, we get
$\Hom_{\Db{B}}(\bar{P}^{0}_{X},B[i])=0$ for all $i> r$ by the above
exact sequence. Therefore,
$$
\End^{i}_{B}(\bar{P}_{X}^{0},B)\simeq
\Hom_{\Db{B}}(\bar{P}^{0}_{X},B[i])=0, \;\;\;\text{for}\;\;\; i>
r.$$ This implies that $\bar{P}_{X}^{0}\in_B\mathcal {X}^{>r}$.
$\square$

Choose an $A$-module $X\in_{A}\mathcal {X}^{>r}$, by Lemma
\ref{2.3}, we know that $F(X)$ is isomorphic in $\Db{B}$ to a
radical complex of the form
$$
0\ra \bar{P}_{X}^{0}\ra \bar{P}_{X}^{1} \ra \cdots\ra
\bar{P}_{X}^{n-1}\ra \bar{P}_{X}^{n}\ra 0
$$
such that $\bar{P}_{X}^{0}\in _{B}\mathcal {X}^{>r}$ and
$\bar{P}_{X}^{i}$ are projective $B$-modules for $1\leq i\leq n$. In
the following, we try to define an additive functor $\underline{F}:
\underline{_{A}\mathcal {X}^{>r}}\ra \underline{_{B}\mathcal
{X}^{>r}}$, where $\underline{_{A}\mathcal {X}^{>r}}$ denotes the
stable category of $_{A}\mathcal {X}^{>r}$ , in which objects are
the same as the objects of $_{A}\mathcal {X}^{>r}$ and, for two
objects $X,Y$ in $\underline{_{A}\mathcal {X}^{>r}}$, their morphism
set is the quotient of $\Hom_{\mathcal {X}_{A}^{>r}}(X,Y)$ modulo
the homomorphisms that factors through projective modules.

\begin{Lem} Let $F:\Db{A}\lra \Db{B}$
be a derived equivalence. Then there is an additive functor
$\underline{F}: \underline{_{A}\mathcal
{X}^{>r}}\ra\underline{_{B}\mathcal {X}^{>r}}$ sending $X$ to
$\bar{P}_{X}^{0}$, such that the following diagram
$$\xymatrix{
      \underline{_{A}\mathcal {X}^{>r}}\ar[r]^(.35){{\rm can}}\ar[d]_{\underline{F}} &
      \Db{A}/\Kb{\pmodcat{A}}\ar[d]^{{F}}\\
      \underline{_{B}\mathcal {X}^{>r}}\ar[r]^(.35){{\rm can}} & \Db{B}/\Kb{\pmodcat{B}}
    }$$
   is commutative up to natural isomorphism.
\end{Lem}
{\bf Proof.} By composing of the embedding functor $_{A}\mathcal
{X}^{>r}\hookrightarrow \Db{A}$ with the localization functor
$\Db{A}\ra \Db{A}/\Kb{\pmodcat{A}}$, we obtain a natural functor
$_{A}\mathcal {X}^{>r}\ra \Db{A}/\Kb{\pmodcat{A}}$. Since the
projective $A$-module is sending to zero in
$\Db{A}/\Kb{\pmodcat{A}}$, we get a canonical functor
$\underline{_{A}\mathcal {X}^{>r}}\ra \Db{A}/\Kb{\pmodcat{A}}$.
There is a functor between $\Db{A}$ and $\Db{B}/\Kb{\pmodcat{B}}$
which is the composition of $F:\Db{A}\lra \Db{B}$ and the
localization functor $\Db{B}\lra\Db{B}/\Kb{\pmodcat{B}}$. Since $F$
is an equivalence, we see that
$F(\Kb{\pmodcat{A}})=\Kb{\pmodcat{B}}$ by Theorem \ref{2.1}. Thus,
there is a functor between $\Db{A}/\Kb{\pmodcat{A}}$and
$\Db{B}/\Kb{\pmodcat{B}}$ induced by $F$, which we also denoted by
$F$. In the following, we will show that the above diagram is
commutative up to natural isomorphism.

For each $f:X\ra Y$ in $_{A}\mathcal {X}^{>r}$, we denote by
$\underline{f}$ the image of $f$ in $\underline{_{A}\mathcal
{X}^{>r}}$. By Lemma \ref{2.3}, we have a distinguished triangle
$$\bar{P}^{+}_{X}\stackrel{i_X}\ra
F(X)\stackrel{j_X}\ra\bar{P}^{0}_{X}\stackrel{m_X}\ra\bar{P}^{+}_{X}[1]\quad
\text{in}\quad \Db{B}.$$ Moreover, for each $f:X\ra Y$ in
$_{A}\mathcal {X}^{>r}$, there is a commutative diagram in $\Db{B}$
$$\xymatrix{
\bar{P}^{+}_{X}\ar^{i_{X}}[r]\ar^{\alpha_{f}}[d]
&F(X)\ar^{j_{X}}[r]\ar^{F(f)}[d]
  & \bar{P}^{0}_{X}\ar^{m_{X}}[r]\ar^{\beta_{f}}[d]&\bar{P}^{+}_{X}[1]
  \ar^{\alpha_{f[1]}}[d]\\
\bar{P}^{+}_{Y}\ar^{i_{Y}}[r] & F(Y)\ar^{j_{Y}}[r]
  & \bar{P}^{0}_{Y}\ar^{m_{Y}}[r]& \bar{P}^{+}_{Y}[1] .}$$
Since $\Hom_{\Db{B}}(\bar{P}^{+}_{X},\bar{P}^{0}_{Y})\simeq
\Hom_{\Kb{B}}(\bar{P}^{+}_{X},\bar{P}^{0}_{Y})=0$, it follows that
$i_{X}F(f)j_{Y}=0$. Then there exists a homomorphism $\alpha_{f}:
\bar{P}^{+}_{X}\ra \bar{P}^{+}_{Y}$. Note that $B$-mod is fully
embedding into $\Db{B}$, hence $\beta_{f}$ is a morphism of
$B$-modules which is in $_{B}\mathcal {X}^{>r}$. If there is another
morphism $\beta'_{f}$ such that $j_{X}\beta'_{f}=F(f)j_{Y}$, then
$j_{X}(\beta_{f}-\beta'_{f})=0$. Thus $\beta_{f}-\beta'_{f}$ factors
through $\bar{P}^{+}_{X}[1]$. There is a distinguished triangle
$$\bar{P}^{1}_{X}[-1]\ra\tau_{\geq
1}(\bar{P}^{+}_{X}[1])\stackrel{a}\ra\bar{P}^{+}_{X}[1]\stackrel{b}\ra\bar{P}^{1}_{X}\quad\text{in}\quad
\Db{A}.$$ Suppose that $\beta_{f}-\beta'_{f}=gh$, where $g: X\ra
\bar{P}^{+}_{X}[1]$ and $h:\bar{P}^{+}_{X}[1]\ra Y$. Since
$\Hom_{\Db{A}}(\tau_{\geq
1}(\bar{P}^{+}_{X}[1]),Y)\simeq\Hom_{\Kb{A}}(\tau_{\geq
1}(\bar{P}^{+}_{X}[1]),Y)=0$, it follows that $ah=0$. Then there is
a map $x: \bar{P}^{1}_{X}\ra Y$, such that $h=bx$. Thus, we get
$\beta_{f}-\beta'_{f}=gbx$, which implies that
$\beta_{f}-\beta'_{f}$ factors through a projective $B$-module.
Therefore, the morphism $\underline{\beta_{f}}$ in
$\Hom_{\underline{\mathcal
{X}_{B}^{>r}}}(\bar{P}^{0}_{X},\bar{P}^{0}_{Y})$is uniquely
determined by $f$.

Let $f: X\ra Y$ and $g: Y\ra Z$ be morphisms in $_{A}\mathcal
{X}^{>r}$. Then there are commutative diagrams as follows:
$$\xymatrix{
\bar{P}^{+}_{X}\ar^{i_{X}}[r]\ar^{\alpha_{fg}}[d]
&F(X)\ar^{j_{X}}[r]\ar^{F(fg)}[d]
  & \bar{P}^{0}_{X}\ar^{m_{X}}[r]\ar^{\beta_{fg}}[d]&\bar{P}^{+}_{X}[1]
  \ar^{\alpha_{fg[1]}}[d]\\
\bar{P}^{+}_{Z}\ar^{i_{Z}}[r] & F(Z)\ar^{j_{Z}}[r]
  & \bar{P}^{0}_{Z}\ar^{m_{Z}}[r]& \bar{P}^{+}_{Z}[1] }$$
  and $$\xymatrix{
\bar{P}^{+}_{X}\ar^{i_{X}}[r]\ar^{\alpha_{f}}[d]
&F(X)\ar^{j_{X}}[r]\ar^{F(f)}[d]
  & \bar{P}^{0}_{X}\ar^{m_{X}}[r]\ar^{\beta_{f}}[d]&\bar{P}^{+}_{X}[1]
  \ar^{\alpha_{f[1]}}[d]\\
\bar{P}^{+}_{Y}\ar^{i_{Y}}[r]\ar^{\alpha_{g}}[d] &
F(Z)\ar^{j_{Y}}[r]\ar^{F(g)}[d]
  & \bar{P}^{0}_{Z}\ar^{m_{Y}}[r]\ar^{\beta_{f}}[d]& \bar{P}^{+}_{Z}[1]\ar^{\alpha_{g[1]}}[d]
  \\\bar{P}^{+}_{Z}\ar^{i_{Z}}[r] & F(Z)\ar^{j_{Z}}[r]
  & \bar{P}^{0}_{Z}\ar^{m_{Y}}[r]& \bar{P}^{+}_{Z}[1] .}$$

Then we have $F(fg)j_{Z}=j_{X}\beta_{fg}$ and
$F(f)F(g)j_{Z}=F(fg)j_{Z}=j_{X}\beta_{f}\beta_{g}$. Therefore,
$j_{X}(\beta_{fg}-\beta_{f}\beta_{g})=0$. By the uniqueness of
$\underline{\beta_{fg}}$, we have
$\underline{\beta_{fg}}=\underline{\beta_{f}}$
$\underline{\beta_{g}}$. Moreover, if $X$ is a projective
$A$-module, then by Lemma \ref{2.3}, we know that $F(X)$ is
isomorphic in $\Db{B}$ to a radical complex of the form
$$
0\ra \bar{P}_{X}^{0}\ra \bar{P}_{X}^{1} \ra \cdots\ra
\bar{P}_{X}^{n-1}\ra \bar{P}_{X}^{n}\ra 0
$$
such that $\bar{P}_{X}^{i}$ are projective $B$-modules for $0\leq
i\leq n$. Thus, if $f$ factors through a projective $A$-module, then
we see that $\beta_{f}$ also factors through a projective
$B$-module.

For each $X\in_{A}\mathcal {X}^{>r}$, we define
$\underline{F}(X):=\bar{P}_{X}^{0}$. Set
$\underline{F}(\underline{f})=\underline{\beta_{f}}$, for each
$\underline{f}\in\Hom_{\underline{\mathcal {X}_{A}}^{>r}}(X,Y)$.
Then $\underline{F}$ is well-defined and an additive functor. The
last statement is discussed in \cite[Proposition 3.5]{P}, we omit it
here. $\square$

The next lemma is useful in our proof of the main result.

\begin{Lem}For $X\in _{A}\mathcal{X}^{>r}$, we have:

For each positive integer $k>r$, there is an isomorphism

$$
\beta_k:
\Hom_{\Db{A}}(X,X[k])\ra\Hom_{\Db{B}}(\underline{F}(X),\underline{F}(X)[k])
$$
Here we denote the image of $g$ under $\beta_k$ by
$\beta_k(g)$.\label{2.4}
\end{Lem}

{\bf Proof.} For $X\in _{A}\mathcal{X}^{>r}$, by Lemma \ref{2.3},
$F(X)=\bar{P}^{\bullet}_{X}$ is isomorphic in $\Db{B}$ to a complex
of the form
$$
0\ra \bar{P}_{X}^{0}\ra \bar{P}_{X}^{1} \ra \cdots\ra
\bar{P}_{X}^{n-1}\ra \bar{P}_{X}^{n}\ra 0
$$
with $\bar{P}_{X}^{0}\in_B\mathcal {X}^{>r}$. Consequently, there is
a distinguished triangle in $\Db{B}$
$$
\bar{P}^{+}_{X}\ra\bar{P}^{\bullet}_{X}\ra\bar{P}^{0}_{X}\ra
\bar{P}^{+}_{X}[1],
$$
where $\bar{P}^{+}_{X}$ is the complex $0\ra\bar{P}_{X}^{1} \ra
\cdots\ra \bar{P}_{X}^{n-1}\ra \bar{P}_{X}^{n}\ra 0$.

For a morphism $f:X\ra X[k]$, and it is easy to see that
$i_{X}F(f)j_{X}[k]\in\Hom_{\Db{B}}(\bar{P}^{+}_{X},\bar{P}^{0}_{X}[k])\simeq
\Hom_{\Kb{B}}(\bar{P}^{+}_{X},\bar{P}^{0}_{X}[k])=0$. Then there is
a map $b_f: \bar{P}^{0}_{X}\ra \bar{P}^{0}_{X}[k]$, we can form the
following commutative diagram
$$
\xymatrix{ \bar{P}^{+}_{X}\ar^{i_{X}}[r]\ar^{a_{f}}[d]
&F(X)\ar^{j_{X}}[r]\ar^{F(f)}[d]
  & \bar{P}^{0}_{X}\ar^{m_{X}}[r]\ar@{.>}^{b_f}[d]&\bar{P}^{+}_{X}[1]
  \ar^{a_{f[1]}}[d]\\
\bar{P}^{+}_{X}[k]\ar^{i_{X}[k]}[r] & F(X)[k]\ar^{j_{X}[k]}[r]
  & \bar{P}^{0}_{X}[k]\ar^{m_{X}}[r]& \bar{P}^{+}_{X}[k+1] .}
$$
We claim that the morphism $b_f$ is uniquely determined by the above
commutative diagram. In fact, if there is another map $b'_f$ such
that $j_{X}b'_f=F(f)j_{X}[k]$, then we get $j_{X}(b_f-b'_f)=0$.
Therefore, $b_f-b'_f$ factors through $\bar{P}^{+}_{X}[1]$. Since
$\Hom_{\Db{B}}(\bar{P}^{+}_{X}[1],\bar{P}^{0}_{X}[k])\simeq
\Hom_{\Kb{B}}(\bar{P}^{+}_{X},\bar{P}^{0}_{X}[k-1])=0$, we have
$b_f-b'_f=0$. Hence, $b_f=b'_f$. Thus, we define a morphism
$$
\beta_k:
\Hom_{\Db{A}}(X,X[k])\ra\Hom_{\Db{B}}(\underline{F}(X),\underline{F}(X)[k]),
$$
by sending $f$ to $b_f$. Next, we will show that $\beta_k$ is an
isomorphism.

Firstly, it is injective. Assume that $\beta_k(f)=b_f=0$. Then
$F(f)j_{X}[k]=0$, and consequently, $F(f)$ factors through
$\bar{P}^{+}_{X}[k]$. It follows that $GF(f)$ factors through
$G(\bar{P}^{+}_{X})[k]$, that is, the map $f: X\ra X[k]$ factors
through $G(\bar{P}^{+}_{X})[k]$, say $f=xy$, for some $x: X\ra
G(\bar{P}^{+}_{X})[k] $ and $y: G(\bar{P}^{+}_{X})[k]\ra X[k]$.
Since $\Hom_{\Db{A}}(\cpx{Q_X},X)\simeq\Hom_{\Kb{A}}(\cpx{Q_X},X)$,
we deduce that $y$ can be chosen to be a chain map. Set
$G(\bar{P}^{+}_{X})=Q^{\bullet}_{X}$. Then $G(\bar{P}^{+}_{X})$ is a
radical projective bounded complex $Q^{\bullet}_{X}$ of the form
$$
0\ra Q^{-m}\cdots\ra Q^{-1}_{X}\ra Q^{0}_{X}\ra Q^{1}_{X}\ra 0,
$$
where $m$ is a positive integer. Indeed, by the distinguished
triangle
$\bar{P}^{+}_{X}\ra\bar{P}^{\bullet}_{X}\ra\bar{P}^{0}_{X}\ra
\bar{P}^{+}_{X}[1],$ we get $H^i(G(\bar{P}^{+}_{X}))=0$ for $i>1$,
where $H^i(G(\bar{P}^{+}_{X}))$ is  $i$-th cohomology group of
$G(\bar{P}^{+}_{X})$. We claim that
$$
\Hom_{\Db{A}}(X,Q^{\bullet}_{X}[k])\simeq\Hom_{\Kb{A}}(X,Q^{\bullet}_{X}[k]).$$
So, it suffices to show that for the complex $Q^{\bullet}_{X}$ of
the form $ 0\ra Q_{X}^{-1}\ra Q_{X}^{0}\ra0$, we get the result.
There is a distinguished triangle $$(\ast)\quad\quad
Q_{X}^{-1}[k]\ra Q_{X}^{0}[k]\ra Q_{X}^{\bullet}[k]\ra
Q_{X}^{-1}[k+1]\quad \text{in} \quad \Kb{A}.$$ Applying the functors
$\Hom_{\Kb{A}}(X,-)$, $\Hom_{\Db{A}}(X,-)$ to ($\ast$), we obtain
the following commutative diagram
$$\xymatrix{
\Hom_{\Kb{A}}(X,Q_{X}^{-1}[k])\ar[r]\ar^{\simeq}[d] &
\Hom_{\Kb{A}}(X,Q_{X}^{0}[k])\ar[r]\ar^{\simeq}[d]
  & \Hom_{\Kb{A}}(X,Q_{X}^{\bullet}[k])\ar[r]\ar[d]& \Hom_{\Kb{A}}(X,Q_{X}^{-1}[k+1]) \ar^{\simeq}[d]\\
 \Hom_{\Db{A}}(X,Q_{X}^{-1}[k])\ar[r] &
\Hom_{\Db{A}}(X,Q_{X}^{0}[k])\ar[r]
  & \Hom_{\Db{A}}(X,Q_{X}^{\bullet}[k])\ar[r]& \Hom_{\Db{A}}(X,Q_{X}^{-1}[k+1]).}$$ Since
$\End^{i}_{A}(X,A)=0$ for $i>r$, it follows that
$\Hom_{\Db{A}}(X,Q_{X}^{-1}[k+1])=0$ for $k>r$. Moreover,
$\Hom_{\Kb{A}}(X,Q_{X}^{-1}[k+1])=0$. We thus get
$\Hom_{\Db{A}}(X,Q^{\bullet}_{X}[k])\simeq\Hom_{\Kb{A}}(X,Q^{\bullet}_{X}[k])$.
Therefore, $x$ is chosen to be a chain map. Consequently, we see
that $f=xy=0$. This shows that $\beta_k$ is injective.

Next, we can prove that $\beta_k$ is surjective. For a map $b:
\bar{P}^{0}_{X}\ra\bar{P}^{0}_{X}[k]$, we have
$j_Xbm_X[k]\in\Hom_{\Db{B}}(F(X),\bar{P}^{+}_{X}[1])$. It follows
from Lemma \ref{PX} that
$$\Hom_{\Db{B}}(F(X),\bar{P}^{+}_{X}[k+1])\simeq\Hom_{\Db{A}}(X,G(\bar{P}^{+}_{X})[k+1])=0 \quad \mbox{for}\quad k>r.$$
Then there is a map $c: F(X)\ra F(X)[k]$ such that
$cj_{X}[k]=j_{X}b$. Since $F$ is an equivalence, it follows that
$c=F(f)$ for some $f: X\ra X[k]$. Hence, $b=\beta_k(f)$. Therefore,
$\beta_k$ is surjective. $\square$

We now have all the ingredients to complete the proof of our main
theorem.

\medskip
{\bf Proof of Theorem \ref{thm1}.} We assume that the generalized
Auslander-Reiten conjecture is true for $B$. If $X$ is an $A$-module
which satisfies $\Ext_A^i(X,X)=0=\Ext_A^i(X,A)$ for $i> r$, then it
follows from Lemma \ref{2.3} that,
$\underline{F}(X)=\bar{P}_{X}^{0}$ satisfies
$\Ext_B^{i}(\bar{P}_{X}^{0},B)=0$ for $i>r$. By Lemma \ref{2.4}, it
is easy to see that, for $i>r$
$$\Ext_B^i(\bar{P}_{X}^{0},\bar{P}_{X}^{0})
\simeq\Hom_{\Db{B}}(\bar{P}_{X}^{0},\bar{P}_{X}^{0}[i])\simeq\Hom_{\Db{A}}(X,X[i])=0.$$
Since we assume that $B$ satisfies the generalized Auslander-Reiten
conjecture, we see that the $\pd(\bar{P}_{X}^{0})\leq r$. We can
take a projective resolution $P^\bullet_{\bar{P}_{X}^{0}}$ of
$\bar{P}_{X}^{0}$. Therefore, by the distinguished triangle
$\bar{P}^{0}_{X}[-1]\ra\bar{P}^{+}_{X}\ra\bar{P}^{\bullet}_{X}\ra\bar{P}^{0}_{X},$
we can take a projective resolution $P^{\bullet}_{\bar{\cpx{P}_X}}$
of $\bar{\cpx{P}_X}$ by the mapping cone of
$P^\bullet_{\bar{P}_{X}^{0}}[-1]$ and $\bar{P}^{+}_{X}$. Thus, we
get $P^{\bullet}_{\bar{\cpx{P}_X}}\in\Kb{\pmodcat{B}}$. It follows
that $X\simeq G(\bar{\cpx{P}_X})\simeq
G(P^{\bullet}_{\bar{\cpx{P}_X}})$. Then it is easy to see that
$\pd_A(X)<\infty$. Let $0\ra P_s\ra\cdots \ra P_r\ra \cdots\ra
P_0\ra X\ra 0$ be a projective resolution of $X$. Since
$\Ext_A^i(X,A)=0$ for $i> r$, we conclude that
$$0\ra \Hom_A(\Omega^r(X),A)\ra\cdots\ra \Hom_A(P_s,A) \ra 0$$
is a split exact sequence. Then, we get $\Hom_A(\Omega^r(X),A)$ is a
projective $A^{op}$-module and consequently, $\Omega^r(X)$ is a
projective $A$-module. It follows that $\pd_A(X)\leq r$.

Similarly, we can prove that the converse is also true.  $\square$

As a corollary of Theorem \ref{thm1}, we re-obtain the following
result of Wei \cite[Theorem 3.7]{W2}.
\begin{Koro} $\rm\cite[Theorem\,3.7]{W2}$ Let $A$ be an Artin algebra and $T$ be a tilting
$A$-module with $\End_A(T)=B$. Then $A$ satisfies the generalized
Auslander-Reiten conjecture if and only if so does $B$.

\end{Koro}

\bigskip

\bigskip


{\footnotesize
\begin{thebibliography}{99}

\bibitem{MR}{{\sc M. Auslander} and {\sc I. Reiten},
On a generalized version of the Nakayama Conjecture, \emph{Proc.
Amer. Math. Soc.} \textbf{52}(1975), 69-74.}

\bibitem{ARS}{{\sc M. Auslander, I. Reiten} and {\sc S.O.Smal\o},
\emph{Representation thoery of Artin algebras}. Cambridge University
Press, 1995.}

\bibitem{CH}{{\sc L. W. Christensen} and {\sc H. Holm},
Algebras that satisfy Auslander's condition on vanishing of
cohomology. \emph{ Math. Z.} \textbf{265} (2010), 21-40.}

\bibitem{Ha1}{{\sc D. Happel},
{\it Triangulated Categories in the Representation Theory of Finite
Dimensional Algebras}. Cambridge University Press, Cambridge. 1988.}

\bibitem{Ha2}{{\sc D. Happel},
Reduction techniques for homological conjecture. \emph{Tsukuba J.
Math.}\textbf{17}(1993),115-130.}

\bibitem{H}{{\sc R. Hartshorne},
{\it Residues and Duality}. Lecture Notes in Math. {\bf20},
Springer-Verlag, 1966.}

\bibitem{hx} {{\sc W. Hu} and {\sc C. C. Xi}, Derived equivalences and
stable equivalences of Morita type, I. \emph{Nagoya Math. J.}
\textbf{200}(2010), 107-152.}

\bibitem{HX3}{{\sc W. Hu} and {\sc C. C. Xi},
Derived equivalences for $\Phi$-Auslander-Yoneda algebras. Preprint,
available at :
http://math.bnu.edu.cn/~ccxi/Papers/Articles/xihu-4.pdf, 2009.}


\bibitem{HL}{{\sc C. Huneke} and {\sc G.J. Leuschke}, On a conjecture of Auslander and
Reiten, \emph{J. Algebra} \textbf{275} (2004), 781-790.}

\bibitem{ke}{{\sc B. Keller}, Derived DG categories. \emph{Ann. Sci.
\'Ecole Norm. Sup.} \textbf{27}(1994), 63-102.}

\bibitem{M}{{\sc B. J. M$\ddot{u}$ller},
The classification of algebras by dominant dimension. \emph{Can. J.
Math.} \textbf{20} (1968), 398-409.}

\bibitem{P}{{\sc S. Y. Pan},
Derived equivalences for Cohen-Macaulay Auslander algebras.
\emph{J.Pure and Appl. Algebra} \textbf{216} (2011), 355-363.}

\bibitem{P1}{{\sc S. Y. Pan},
Derived equivalences for $\Phi$-Cohen-Macaulay Auslander-Yoneda
algebras. Submitted.}

\bibitem{PX}{{\sc S. Y. Pan} and {\sc C. C. Xi},
Finiteness of finitistic dimension is invariant under derived
equivalences. {\it J. Algebra} \textbf{322} (2009), 21-24.}

\bibitem{Ri1}{{\sc J. Rickard}, Morita theory for derived categories.
\emph{J. London Math. Soc}. \textbf{39} (1989) 436-456.}

\bibitem{R2}{{\sc J. Rickard}, Derived categories and stable equivalences. \emph{J. Pure Appl.
Algebra} \textbf{61}(1989)303-317.}

\bibitem{Ri3}{{\sc J. Rickard},
Derived equivalences as derived functors. \emph{J. London Math.
Soc.} \textbf{43}(1991), 37-48.}

\bibitem{T}{{\sc T. Araya}, Auslander-Reiten conjecture on Gorenstein rings.
\emph{Proc. Amer. Math. Soc.} \textbf{137} (2010), 1941-1944.}

\bibitem{W1}{{\sc J. Q. Wei}, A note on Auslander bounds.
Preprint, available at :http://arxiv.org/abs/0802.1085v2.}

\bibitem{W2}{{\sc J. Q. Wei},
Generalized Auslander-Reiten conjecture and tilting equivalences.
\emph{Proc. Amer. Math. Soc.} \textbf{138} (2010), 1581-1585.}

\bibitem{W3}{{\sc J. Q. Wei},
Tilting invariance of the Auslander-Reiten conjecture. \emph{Math.
Res. Lett. } \textbf{17} (2010), 171-176.}

\bibitem{C}{{\sc C. C. Xi},
On the finitistic dimension conjecture, I: related to
representation-finite algebras. \emph{J. Pure Appl.
Alg.}\textbf{193} (2004), 287-305.}
\bibitem{Ver}{{\sc J. Verdier},
Cat\'{e}gories d\'{e}riv\'{e}es, \'{e}tat 0. \, \emph{Lecture Notes
in Math. 569}(1977), Springer, Berlin, 262-311.}

\end{thebibliography} }
\end{document}